\def\ove#1{\overline{#1}}    %
\def\ovs#1#2{\overset{#1}\to{#2}}
     \def\({\left(}                  
     \def\){\right)}          
     \def\[{\left[}       \def\la{\lambda}
     \def\]{\right]}      \def\ffi{\varphi}

                                      \def\ot{\otimes}
     \def\<{\langle}                 \def\wh{\widehat}
     \def\>{\rangle}                 
                 \def\sbs{\subset}
\def\tr{\operatorname{trace}\,}

\def\N{\operatorname{N}}

\def\id{\operatorname{id}}

\def\R{\operatorname{R}}

                     \def\sbs{\subset}

    \def\QQ{$\quad\square$}

\documentclass[12pt,oneside]{amsproc}
\usepackage[T2A]{fontenc}
\usepackage[english]{babel}
\usepackage{amssymb}
\usepackage{amsthm}

\title{Grothendieck-Lidski\v{\i} theorem for subspaces and  factor spaces of  $L_p$-spaces }
\author{Oleg Reinov}
\address{Department of Mathematics and Mechanics, St. Petersburg State University,
Saint Petersburg, RUSSIA.\newline
\phantom{Ao} Abdus Salam School of Mathematical Sciences, 68-B, New Muslim Town, Lahore 54600, PAKISTAN.
}
\email{orein51@mail.ru}
\author{Qaisar Latif}
\address{Abdus Salam School of Mathematical Sciences, 68-B, New Muslim Town, Lahore 54600, PAKISTAN.}

\thanks{%${ }^\maltese$
The research was supported by the Higher  Education Commission of Pakistan.}

\thanks{%${ }^\maltese$
}
\thanks{${ }$ Key words: $s$-nuclear operators, eigenvalue distributions. }

\begin{document}

                          $$ {} $$
\vphantom{} \maketitle

\begin{abstract}
In 1955, A. Grothendieck  has shown that if the linear operator $T$
in a Banach subspace of an $L_\infty$-space is $2/3$-nuclear then
the trace of $T$ is well defined and is equal to the sum of all  eigenvalues $\{\mu_k(T)\}$
of $T.$  \
V.B. Lidski\v{\i} , in 1959, proved his famous theorem on the coincidence of the trace of the $S_1$-operator
in $L_2(\nu)$ with its spectral trace $\sum_{k=1}^\infty \mu_k(T).$            \
We show that {\it
for $p\in[1,\infty]$ and $s\in (0,1]$ with $1/s=1+|1/2-1/p|,$
and for every $s$-nuclear operator $T$ in every subspace of any $L_p(\nu)$-space
the trace of $T$ is well defined and equals the sum of all eigenvalues of $T.$
 }                                                                                                            \
Note that for $p=2$ one has $s=1,$ and  for $p=\infty$ one has $s=2/3.$
 \end{abstract}

\vskip 0.75cm

 In 1955, A. Grothendieck [1] has shown that if the linear operator $T$
in a Banach space  is $2/3$-nuclear then
the trace of $T$ is well defined and is equal to the sum of all  eigenvalues $\{\mu_k(T)\}$
of $T.$  \,
V.B. Lidski\v{\i} [2], in 1959, proved his famous theorem on the coincidence of the trace of the $S_1$-operator
in an (infinite dimensional) Hilbert space
 with its spectral trace $\sum_{k=1}^\infty \mu_k(T).$
  Any Banach space is a subspace of an $L_\infty(\nu)$-space, as well as any Hilbert space
  is a (subspace of) $L_2(\nu)$-space. Also, any Banach space is a factor space of an $L_1(\nu)$-space,
  as well as any Hilbert space is a (factor space of) $L_2(\nu)$-space.
  We obtain the following generalization of these theorems:
  {\it
for $p\in[1,\infty]$ and $s\in (0,1]$ with $1/s=1+|1/2-1/p|,$
and for every $s$-nuclear operator $T$ in every subspace of any $L_p(\nu)$-space
the trace of $T$ is well defined and equals the sum of all eigenvalues of $T.$}
 Note that for $p=2$ one has $s=1,$ and  for $p=\infty$ one has $s=2/3.$

  \vskip 0.33cm

\centerline{\bf \S1. Definitions and a theorem}

 \vskip 0.23cm

All the terminology and facts (now classical), given here without any explanations, can be found
in [7--10].
         \vskip 0.02cm
  Let $X,Y$ be Banach spaces. For $s\in (0,1],$ denote by $X^*\wh\ot_s Y$
  the completion of the tensor product $X^*\ot Y$  (considered as a linear space
  of all finite rank operators) with respect to the quasi-norm
  $$
  ||z||_s:= \inf \{\(\sum_{k=1}^N ||x'_k||^s\, ||y_k||^s\)^{1/s}:\ z=\sum_{k=1}^N x'_k\ot y_k\}.
  $$
  Let $\Phi_p,$ for $p\in[1,\infty],$ be the ideal of all operators which can be factored through a subspace af an $L_p$-space.
  Put $N_s(X,Y):= $ image of $X^*\wh\ot_s Y$ in the space $L(X,Y)$ of all bounded linear transformations under the canonical factor map
  $X^*\wh\ot_s Y\to N_s(X,Y)\sbs L(X,Y).$ We consider the (Grothendieck) space $N_s(X,Y)$  of all
  $s$-nuclear operators from $X$ to $Y$ with the natural quasi-norm, induced from $X^*\wh\ot_s Y.$

  Finally, let $\Phi_{p,s}$ (respectively, $\Phi_{s,p})$ be the quasi-normed product $N_s\circ \Phi_p$
  (respectively, $\Phi_p\circ N_s)$ of the corresponding ideals equipped with the natural quasi-norm $\nu_{p,s}$
  (respectively, $\nu_{s,p})$: if $A\in N_s\circ \Phi_p(X,Y)$ then
  $A=\ffi\circ T$ with  $T=\beta\alpha\in \Phi_p,$ $\ffi=\delta\Delta\gamma\in N_s$ and
  $$
  A: X \overset\alpha\to X_p\overset\beta\to Z \overset\gamma\to c_0\overset\Delta\to l_1\overset\delta\to Y,
  $$
  where all maps are continuous and linear, $X_p$ is a subspace of an $L_p$-space, constructed on a measure space, and
  $\Delta$ is a diagonal operator with the diagonal from $l_s.$
  Thus, $A=\delta \Delta\gamma\beta\alpha$ and $A\in N_s.$
  Therefore, if $X=Y,$    the spectrum of $A,$\, $sp\,(A),$
  is at most countable with only possible limit point zero. Moreover, $A$ is a Riesz operator with eigenvalues
  of finite algebraic multiplicities and $sp\, (A)\equiv sp\, (B),$ where $B:= \alpha\delta\Delta\gamma\beta: X_p\to X_p$
 is an $s$-nuclear operator, acting in a subspace of an $L_p$-space.

  \vskip 0.1cm

  Let $T$ be an operator between Banach spaces  $Y$ and $W.$ % $z\in X^*\wh\ot Y.$
 The operator ${\mathbf1}\ot T: X^*\wh\ot_s Y\to X^*\wh\ot_s W$ is well defined and can be considered also as
 an operator from $X^*\wh\ot_s Y$ into $X^*\wh\ot W$ (the Grothendieck projective tensor product),
 the last space having the space $L(W,X^{**})$ as dual.

\vskip 0.1cm

{\bf Definition.}\
We say that $T$ {\it possesses the property $AP_s$}\,
(written down as "$T\in AP_s$") if  for every $X$ and any tensor element $z\in X^*\wh\ot_s Y$
the operator $T\circ z: X\to W$ is zero iff the corresponding tensor $({\mathbf1}\ot T)(z)$
is zero as an element of the space $X^*\wh\ot W.$   If $Y=W$ and $T$ is the identity map, we write just
$Y\in AP_s$ (the approximation property of order $s).$
\vskip 0.1cm

This is equivalent to the fact that if $z\in X^*\wh\ot_s Y$  then it follows from
$$
 \tr (\mathbf1\ot T)(z)\circ R=0, \ \ \forall\, R\in W^*\ot X
$$
that $\tr U\circ (\mathbf1\ot T)(z)=0$ for every  $U\in L(W,X^{**}).$
There is a simple characterization of the condition $T\in AP_s$
in terms of the approximation of $T$ on some sequences of the space $Y,$
but we omit it now, till the next time.
We need here only one example which is crucial for our note (other
examples, as well as more general applications will appear elsewhere).
     \vskip 0.1cm

     {\bf Example.}\
     Let $s\in (0,1],$ $p\in [1,\infty]$ and $1/s=1+|1/p-1/2|.$
     Any subspace as well as any factor space
      of any $L_p$-space have the property $AP_s$
     (this means that, for that space $Y,$  $\id_Y\in AP_s).$
     Thus, in the case of such a space $Y,$ we have the quasi-Banach equality
     $X^*\wh\ot_s Y=N_s(X,Y),$ whichever the space $X$ was.

       \vskip 0.1cm

{\bf Lemma.}\
Let $s\in (0,1],$ $p\in [1,\infty]$ and $1/s=1+|1/2-1/p|.$
Then the system of all eigenvalues (with their algebraic multiplicities)
of any operator $T\in N_s(Y,Y),$ acting in any subspace $Y$ of any
$L_p$-space, belongs to the space $l_1.$ The same is true for the factor spaces
of $L_p$-spaces.

           \vskip 0.1cm

           {\bf Corrolary.}\
           If    $s\in (0,1],$ $p\in [1,\infty]$ with $1/s=1+|1/2-1/p|$
           then the quasi-normed ideals $\Phi_{p,s}$ and $\Phi_{s,p}$
           are of (spectral) type $l_1.$
   \vskip 0.1cm

  {\bf Theorem.}\
Let $Y$ be a subspace or a factor space of  an $L_p$-space,
$1\le p\le \infty.$ If $T\in N_s(Y,Y),$\,
$1/s=1+|1/2-1/p|,$   \,
then

1.\, the (nuclear) trace  of $T$ is well defined,

2.\, $\sum_{n=1}^\infty |\la_n(T)|<\infty,$ where
$\{\la_n(T)\}$ is the system of all eigenvalues of the operator $T$
(written in according to their algebraic multiplicities)

and
$$
 \tr T= \sum_{n=1}^\infty \la_n(T).
$$

\vskip 0.3cm
%\newpage

\centerline{\bf \S2. Proofs}

 \vskip 0.2cm

{\it Proof}\ of Lemma.    \
Let $Y$ be a subspace or a factor space of an $L_p$-space
and $T\in N_s(Y,Y)$ with an s-nuclear representation
$$
 T=\sum_{k=1}^\infty \mu_k y'_k\ot y_k,
$$
where $||y'_k||, ||y_k||=1$ and $\mu_k\ge 0,$  $\sum_{k=1}^\infty \mu_k^s<\infty.$
%Let-... omited
The operator $T$ can be factored in the following way:
$$
 T: Y\overset{A}\longrightarrow l_\infty \overset{\Delta_{1-s}}\longrightarrow l_r
 \overset{j}\hookrightarrow c_0 \overset{\Delta_s}\longrightarrow l_1\overset{B}\longrightarrow Y,
$$
where $A$ and $B$ are linear bounded, $j$ is the natural injection, $\Delta_s\sim(\mu_k^s)_k$ and
$\Delta_{1-s}\sim(\mu_k^{1-s})$ are the natural diagonal operators from $c_0$ into $l_1$ and
from $l_\infty$ into $l_r,$ respectively. Here, $r$ is defined via the conditions
 $1/s=1+|1/p-1/2|$ and $\sum_k \mu_k^s<\infty:$
 we have to have $\sum_k \mu_k^{(1-s)r}<\infty,$ for which $(1-s)r=s$ is good. Therefore, put
 $1/r=1/s-1,$ or $1/r=|1/p-1/2|.$

        from now onward in the proof we assume (surely, without loss of generality) that
        $p\ge2.$
Then $1/r=1/2-1/p$ and $r(1-s)=s.$
Note that if $s=1$ then $r=\infty, p=2$ and $j\Delta_{1-s}\equiv j;$
and if $s=2/3$ then  $r=2, p=\infty$ and $\Delta_{1-s}\sim (\mu_k^{1/3})_k\in l_2.$

Now, let us factorize the diagonal $\Delta_s$ as
$\Delta_s=\Delta_2\Delta_1:\, c_0 \overset{\Delta_{1}}\longrightarrow l_2  \overset{\Delta_{2}}\longrightarrow l_1$
in such a (clear) way that diagonals $\Delta_1$ is in $\Pi_2$ and $\Delta^*_2$ is in $\Pi_2$ too, respectively.

{\it Case }\, (i).\ $Y$ is a subspace of an $L_p$-space.
Denoting by $l: Y\hookrightarrow L_p$ an isomorphic embedding of $Y$ into a corresponding $L_p=L_p(\nu),$
we obtain that the map
$\Delta^*_2B^*l^*: \, L_{p'}\overset{l^*}\longrightarrow Y^*\overset{B^*}\longrightarrow l_\infty \overset{\Delta^*_{2}}\longrightarrow l_2$
is of type $\Pi_2,$    so is in $\Pi_p.$
Thus its preadjoint $lB\Delta_2:\, l_2\overset{\Delta_2}\longrightarrow l_1\overset{B}\longrightarrow Y \overset{l}\longrightarrow L_p$
is order bounded and, therefore, $p$-absolutely summing.

{\it Case }\, (ii).\ $Y$ is a factor space of an $L_p$-space.
Denoting by $q: L_p\to Y$ a factor map from a corresponding $L_p=L_p(\nu)$ onto $Y$ and
taking a lifting $Q: l_1\to L_p$ for $B$ with $B=qQ,$
we obtain that the map
$\Delta^*_2Q^*: \, L_{p'}\overset{Q*}\longrightarrow l_\infty \overset{\Delta^*_{2}}\longrightarrow l_2$
is of type $\Pi_2,$    so is in $\Pi_p.$
Thus its pre-adjoint $Q\Delta_2:\, l_2\overset{\Delta_2}\longrightarrow l_1\overset{Q}\longrightarrow L_p$
is order bounded and, therefore, $p$-absolutely summing.
Hence, $B\Delta_2: \, l_2\overset{\Delta_2}\longrightarrow l_1\overset{Q}\longrightarrow L_p\overset{q}\longrightarrow Y$
is also $p$-absolutely summing.

It folows from all that's said that in all the cases our operator $T:Y\to Y$ can be written as a composition:
 $$
  T=U_1U_2U_3\ \text{ with } \ U_3\in\Pi_r, U_2\in \Pi_2, U_1\in \Pi_p,
 $$
all the exponents being not less than 2.
 Now, $1/r+1/2+1/p=(1/2-1/p)+1/2+1/p=1.$
 %Therefore, Lemma is proved.
  \QQ

    \vskip 0.1cm

 {\it Proof}\ of the statement of Example.
  It follows from:

 $(\alpha)$\ every finite dimensional subspace $E$ of any factor space of any $L_p$-space
is $c_p\,(\dim{E})^{|1/2-1/p|}$-complemented.

\noindent
 For more general statements on $AP_s$ and their proofs,
we refer to [4] and [5]; see also an old paper of O.I. Reinov [3] for the idea to apply the projections
in the questions which are under consideration in this note.     \QQ

          \vskip 0.1cm

 {\it Proof}\ of Corrolary.
 Apply Lemma.        \QQ

      \vskip 0.1cm

 {\it Proof}\ of Theorem.
 Apply Lemma, Example, Corrolary and the main  result of M.C. White [6].

     \vskip 0.1cm

 {\it Remark}:\,
 Since finite rank operators are dense in $N_s,$
 Theorem can be proved without referring to the paper of M.C. White;
  but this would take a little bit longer explanations.

   \newpage

%%%%%%%%%%       ++++++++++++++++++++++++++++++++++++++++++++++++++++++++++

%\bigskip
%\bigskip
%\medskip

%%%%%%%%%%%%%%%%%%%%%%%%%%%%%%%%%%%%%%%%%%%%%%%%%%%%%%%%%%%%%
\end{document}